\documentclass[runningheads,a4paper]{llncs}
\usepackage{amssymb}
\usepackage{amsmath}
\usepackage{graphicx}
\usepackage[ruled,nothing]{algorithm}
\usepackage{url}
\usepackage{nicefrac}

\newcommand{\p}{^{\prime}}

\title{On connection between division sequences and presentations of a free group}

\mainmatter              
\institute{CS department, Ben-Gurion University, Beer-Sheva, Israel\\
\email{orlovn@cs.bgu.ac.il}
}
\author{Natalia Vanetik}

\begin{document}
\pagestyle{plain}

\maketitle
\begin{abstract}
  This paper describes a simple method for estimating lower bounds on the number of
  classes of equivalence for a special kind of integer sequences, called
  division sequences. The method is based on adding group structure
  to classes of equivalence and
  studying properties of resulting groups as presentations of free group.
\end{abstract}

\section{Division sequences}
Let $M\in \{\mathbb{Z}_{\{\neq 0\}}, \mathbb{Z}_{\{>0\}}\}$ be a
domain of mapping 
$\mathcal{C}_{p,q}^{M}:M\rightarrow M$ defined as follows.
\begin{equation}\mathcal{C}_{p,q}^{M}(c)=\left \{ \begin{array}{ll} pc+1, &
      \mbox{\sl $q$ does not divide $c$} \\  c/q,
      & \mathrm{otherwise} \end{array}\right .\end{equation}
We call the mapping $\mathcal{C}_{p,q}^{M}$ \textit{division sequence} for parameters $p$ and $q$ and domain $M$.

Division sequence $\mathcal{C}_{p,q}^{M}$ indices following equivalence
relation on $M$: for $c,d\in M$ we say that $c$ is
\textit{equivalent} to $d$, denoted $c\sim d$, if
$d=\mathcal{C}_{p,q}^{M}(c)$; symmetry, reflexivity and transitivity are
added to $\sim$ in order for it to be a proper equivalence relation.

\section{Presentations of free group}
Let us observe monoid $M\in \{\mathbb{Z}_{\{\neq 0,\cdot\}},
\mathbb{Z}_{\{>0,\cdot\}}\}$ and corresponding Grothendieck's
group $F:=K(M)\in \{\mathbb{Q}_{\{\neq 0,\cdot\}}, \mathbb{Q}_{\{>0,\cdot\}}\}$.

We are interested in free group $F$ with countably many generators (primes) and
defining relations that (a) include abelinizing relation and (b) reflect equivalence
relation $\sim$ defined on $M$ by $\mathcal{C}_{p,q}^{M}$.

For fixed $p,q\in M$ we observe following relations:
\begin{align}\label{def-rel-1}&\mbox{\sl $q=1$,}\\
\label{def-rel-2}&P\sim Q \Rightarrow
  P=Q
\end{align}

For given $p$ and $q$, we denote relations generated by \eqref{def-rel-2} by
$R_{p,q}$. Abelinizing relation that we add to $F$ is denoted $XY=YX$ and it
implies that $xy=yx$ for all pairs of generators $x,y$ of $F$. Note that
$R_{p,q}$ is a countable set.

Let $$H_{p,q}:=\langle F\: |\: XY=YX,\: q=1, \: R_{p,q}\:\rangle$$ be  presentation of $F$ (we write $F$ instead of
$F$'s generators for simplicity).
We denote by
$\phi_{p,q}$ the corresponding free group homomorphism defining
$H_{p,q}$. If the choice of $F$ needs to be clarified for $H_{p,q}$, we
write $H_{p,q}^{M}$ for domain $M$ and free group $F=K(M)$.
\begin{property}\label{q^n}$H_{p,q}$ is isomorphic to a quotient
  group of $H_{p,q^{n}}$ for all $n\in
    \mathbb{N}$.\end{property}
\proof By Tietze transformations rules (see, e.g., \cite{tietze}), implication $q=1 \Rightarrow
q^{n}=1$ ensures that $H_{p,q}$ is isomorphic to a quotient group of
$H_{p,q^{n}}$. Since natural homomorphism $\chi:H_{p,q}\rightarrow
H_{p,q^{n}}$, where $\chi(a)=a$, is surjective we have $\mathrm{Im\: }\chi \unlhd
H_{p,q^{n}}$. 
\qed

When we are interested in parameters $p,q$ for which $H_{p,q}\cong \{1\}$,
i.e. is a trivial group, following properties are obvious.
\begin{align}
    &H_{p,q}\cong \{1\} \Longrightarrow
    H_{p,q^{n}}\cong \{1\}\: \mathrm{and}\\
    &H_{p,q}\ncong \{1\} \Longrightarrow H_{p,q^{n}}\ncong
    \{1\}\: \forall n\in \mathbb{N}\nonumber \end{align}
Another indication of structure of $H_{p,q}$ are groups
$H_{pp\p,q}$, where $p\p\in F$.
\begin{property}\label{adding-p}Let $p\p\in \mathrm{Ker\: }\phi_{p,q}$. Then $H_{p,q}$ is
  isomorphic to a quotient group of $H_{pp\p,q}$.\end{property}
\proof Let $pp\p P+1=Q$ imply defining relation $P=Q \in R_{pp\p,q}$.
Since $H_{p,q}$ and $H_{pp\p,q}$ share domain, $pp\p P+1=Q$ also implies
that defining relation $p\p P=Q \in R_{p,q}$. Relations $p\p P=Q, p\p=1 \in R_{p,q}$ 
give us defining relation $P=Q \in R_{p,q}$ by Tietze transformations
rules.
Therefore, each defining relation of $H_{pp\p,q}$ is also a defining
relation of $H_{p,q}$, and we have the property. \qed

Tietze transformations allow us to make a stronger claim under following
restrictions.
\begin{property}\label{adding-q}Let $q\p\in \mathrm{Ker\: }\phi_{p,q} \cap \mathrm{Ker\: }\phi_{p,qq\p}$. Then \\
 \hspace*{5mm}(a) $H_{p,q}$ is isomorphic to a quotient group of $H_{p,qq\p}$,\\
 \hspace*{5mm}(b) $H_{p,qq\p}$ is isomorphic to a quotient group of $H_{p,q}$. 
  \end{property}
\proof Condition $q\p\in \mathrm{Ker\: }\phi_{p,q} \cap \mathrm{Ker\: }\phi_{p,qq\p}$ ensures that
\begin{equation}\label{eq1}q=1,\: q\p=1,\: qq\p=1 \in R_{p,q}\cap R_{p,qq\p}\end{equation}
Let us assume that equality
\begin{equation}\label{eq2}pP+1=q^{n}(q\p)^{m}Q\end{equation}
holds in $\mathbb{Q}$. Then $pP+1=(q\p)^{m}Q$ holds in
$\nicefrac{\mathbb{Q}}{\langle q=1\rangle}$ and thus
\begin{equation}\label{eq3}P=(q\p)^{m}Q\in R_{p,q}\Rightarrow
P=Q\in R_{p,q}\end{equation}
On other hand, equality \eqref{eq2} gives us defining relations
\begin{equation}P=q^{n-\min(m,n)}(q\p)^{m-\min(m,n)}Q\in R_{p,qq\p}
  \Rightarrow P=Q\in R_{p,qq\p}\end{equation} by Tietze
transformations rules because
$q\p=1,q=1\in R_{p,qq\p}$.
Therefore, $H_{p,q}$ is isomorphic to a quotient group of
$H_{p,qq\p}$ and vice versa. \qed
\begin{corollary}\label{adding-q-corollary}If one of $H_{p,q}$, $H_{p,qq\p}$ is finite and Property
  \ref{adding-q} holds for $q\p$, $H_{p,q}\cong H_{p,qq\p}$ by cardinality argument. In particular, $H_{p,q}\cong \{1\}$
    if and only if $H_{p,qq\p}\cong \{1\}$.\qed\end{corollary}
Since moving from group $\mathbb{Q}_{\{\neq 0,\cdot\}}$ to
group $\mathbb{Q}_{\{>0,\cdot\}}$
in terms of defining relations can be done by adding
defining relation $-1=1$, we also have
\begin{equation}\mbox{\sl $H_{p,q}^{\mathbb{Q}_{\{>0,\cdot\}}}$ is
    isomorphic to a quotient group of
  $H_{p,q}^{\mathbb{Q}_{\{\neq 0,\cdot\}}}$}\end{equation}
Let us now observe group $\mathrm{Ker\: }\phi_{p}$ which is a
subgroups of $F$. Since $\mathrm{Ker\: }\phi_{p,q}\subseteq \mathrm{Ker\:
}\phi_{p,q^{n}}$ by Property  \ref{q^n},
we have normal subgroup relation
\begin{equation}\label{inverse-sequence-1}\mathrm{Ker\: }\phi_{p,q} \unlhd \mathrm{Ker\: }\phi_{p,q^{n}}\end{equation}
Then following corollaries are implied by Properties \ref{adding-p}-\ref{adding-q}.
\begin{corollary}\label{inverse-sequence-2}Let  
  $p\p\in \mathrm{Ker\: }\phi_{p,q}$. Then $\mathrm{Ker\: }\phi_{pp\p,q}
  \unlhd \mathrm{Ker\: }\phi_{p,q}$. \qed \end{corollary}
\begin{corollary}\label{inverse-sequence-3}Let $q\p\in M$ 
  such that
  $q\p\in \mathrm{Ker\: }\phi_{p,q}\cap \mathrm{Ker\: }\phi_{p,qq\p}$.
  Then $\mathrm{Ker\: }\phi_{p,q}=\mathrm{Ker\:
  }\phi_{p,qq\p}$ and $\mathrm{Ker\: }\phi_{p,q}\cong \mathrm{Ker\:
  }\phi_{p,qq\p}$. \qed \end{corollary}

\vspace*{2mm}
Structure of $\mathrm{Ker\: }\phi_{p,q}$, if nontrivial, can be partially
disclosed as follows.
Let us fix $p,q\in M$ and define
\begin{equation}\overline{R}_{p,q}:= P\nsim 1 \Rightarrow P=1\end{equation}
to be set of defining relations, where $\sim$ is equivalence relation of
a division sequence $\mathcal{C}_{p,q}^{M}$.

Now we can define group $\overline{H}_{p,q}$ as representation of free
group $F$ as follows.
\begin{equation}\label{overline-G}\overline{H}_{p,q}:=\langle F\: |\:
  XY=YX,\: q=1,\: \overline{R}_{p,q} \rangle\end{equation}
We denote the corresponding homomorphism of $F$ by $\overline{\phi}_{p,q}$.
Since $\mathrm{Im\: }\overline{\phi}_{p,q}\subseteq \mathrm{Ker\:
}\phi_{p,q}$,
the following holds.
\begin{property}$\overline{H}_{p,q}$ is isomorphic to a quotient group of $\mathrm{Ker\:
}\phi_{p,q}$. \qed\end{property}
\begin{property}\label{not-H-subgroup} Let $p\p,q\p\in M$. 
   Then $\nicefrac{\overline{H}_{p,q}}{\langle p\p=1
    \rangle}$ is isomorphic to a quotient group of $\overline{H}_{pp\p,q}$ and
  $\nicefrac{\overline{H}_{p,q}}{\langle q\p=1 \rangle}$ is
  isomorphic to a quotient group of $\overline{H}_{p,qq\p}$.\end{property}
\proof Equivalence relation $P\sim Q$ implied by equality $pp\p P+1=Q$ generates defining relation $P=Q\in
R_{pp\p,q}$, defining relation $p\p P=Q\in R_{p,q}$ in $H_{p,q}$ and
defining relation $P=Q$ of $\nicefrac{H_{p,q}}{\langle p\p=1
  \rangle}$.
Therefore relations $P\sim Q\in \overline{R}_{pp\p,q}$, where $P,Q\nsim 1$,
generate defining relation $P=Q\in R_{p,q}$ and therefore are 
defining relations of $\nicefrac{\overline{H}_{p,q}}{\langle p\p=1 \rangle}$
as well.
  
Similarly, equivalence relation $P\sim Q$, implied by $pp\p P+1=Q$, generates relation
$P=Q\in R_{p,qq\p}$ and relation $P=(q\p)^{i}Q\in R_{p,q}$ for some $i\in
\mathbb{N}\cup \{0\}$. Then 
$P=Q$ is a defining relation of $\nicefrac{H_{p,q}}{\langle q\p=1 \rangle}$
as well. 
Therefore relations $P\sim Q$, where $P,Q\nsim 1$, that
generate $\overline{R}_{p,qq\p}$, also generate
defining relations of $\nicefrac{\overline{H}_{p,q}}{\langle q\p=1 \rangle}$.

Thus, group $\nicefrac{\overline{H}_{p,q}}{\langle p\p=1 \rangle}$ is
isomorphic to a quotient group of $\overline{H}_{pp\p,q}$
and $\nicefrac{\overline{H}_{p,q}}{\langle q\p=1 \rangle}$ is
  isomorphic to a quotient group of $\overline{H}_{p,qq\p}$. \qed
\begin{property}\label{qp-factorization}Let $q\p\in \mathrm{Ker\: }\phi_{p,q} \cap \mathrm{Ker\:
  }\phi_{p,qq\p}$. Then $\nicefrac{\overline{H}_{p,q}}{\langle q\p=1 \rangle}\cong
\nicefrac{\overline{H}_{p,qq\p}}{\langle q\p=1 \rangle}$. \end{property}
\proof By Property \ref{adding-q}, defining relations of groups
$\nicefrac{\overline{H}_{p,q}}{\langle q\p=1 \rangle}$
and $\nicefrac{\overline{H}_{p,qq\p}}{\langle q\p=1 \rangle}$ are generated
by the same set of relations.
Then $P=1$ is a defining relation of $\nicefrac{\overline{H}_{p,q}}{\langle q\p=1 \rangle}$ if
and only if it is a defining relation of
$\nicefrac{\overline{H}_{p,qq\p}}{\langle q\p=1 \rangle}$.
Therefore, $\nicefrac{\overline{H}_{p,q}}{\langle q\p=1 \rangle}\cong
\nicefrac{\overline{H}_{p,qq\p}}{\langle q\p=1 \rangle}$, since
$q=1,q\p=1,qq\p=1$ are defining relations of both groups 
\qed

\section{Equivalence classes of division sequences}
Let us get back to division sequences $\mathcal{C}_{p,q}^{M}$ over domain
$M\in \{\mathbb{Z}_{\{\neq 0\}}, \mathbb{Z}_{\{>0\}}\}$.
Free group presentations $H_{p,q}$ are related to equivalence classes
defined by $\mathcal{C}_{p,q}^{M}$ as follows.
Relation \eqref{def-rel-1} restricts equivalence relation $\sim$ to positive
integer numbers co-prime to $q$, and relation  \eqref{def-rel-2} ensures
that for every $P\in M$ such that $P\sim 1$ we have $P=1$
in $H_{p,q}$ (the inverse implication may not be true). Therefore,
elements of $H_{p,q}$ that are not identity contain
positive integer numbers that are not equivalent to $1$ in the sense of $\sim$.
\begin{corollary}\label{G-1}Let $c,d\in M$. If $c$ and $d$ lie in different
  elements of $H_{p,q}$, then $c\nsim d$. Therefore, the number of
  different equivalence classes of $\mathcal{C}_{p,q}^{M}$ is at least the order of $H_{p,q}$. \qed\end{corollary} %

The structure of groups $\mathrm{Ker\: }\phi_{p,q}$
is somewhat an indication of the structure of integers equivalent
to $1$ in $\mathcal{C}_{p,q}^{M}$.
\begin{corollary}\label{kernel-group}Members of equivalence classes of $\sim$ that lie
  in identity element of $H_{p,q}$ form a group under multiplication. \qed\end{corollary} 
\begin{corollary}\label{not-G-2}If $\mathrm{Ker\: }\phi_{p,q}= F$, then $\mathrm{Ker\: }\phi_{p,q^{n}}= F$ for all
$n\in \mathbb{N}$. \qed\end{corollary}
Let us give several examples of how the above properties can be applied
to specific division sequences.
\begin{example}If $\mathcal{C}_{p,q}^{\mathbb{Z}_{\{>0\}}}$ has just one
  equivalence class, then $H_{p,q}^{\mathbb{Q}_{\{\neq 0,\cdot\}}}$ has order $\leq 3$.\end{example}
\proof Let us assume that $H_{p,q}^{\mathbb{Q}_{\{\neq 0,\cdot\}}}$ has order $\geq 3$. Then
there exist cosets 
$a\mathrm{Ker\: }\phi_{p,q}$, $b\mathrm{Ker\: }\phi_{p,q}$ and
$c\mathrm{Ker\: }\phi_{p,q}$ where $a,c,b\in \mathbb{Q}_{\{<0\}}$.
In this case $ab,ac\in \mathbb{Q}_{\{>0\}}$ imply that $a\mathrm{Ker\: }\phi_{p,q}\cdot b\mathrm{Ker\: }\phi_{p,q}=ab
{Ker\: }\phi_{p,q}={Ker\: }\phi_{p,q}$
and similarly $a\mathrm{Ker\: }\phi_{p,q}\cdot c\mathrm{Ker\: }\phi_{p,q}={Ker\:
}\phi_{p,q}$ - a contradiction. \qed
\begin{example}If $\mathcal{C}_{p,q^{n}}^{M}$, $n\in \mathbb{N}$, 
  has just one equivalence class, then 
  $H_{p,q}^{M}\cong \{1\}$. \qed \end{example}
\begin{example}If $\mathcal{C}_{pp\p,q}^{M}$ contains a single equivalence
  class, then $H_{p,q}^{M}\cong \{1\}$. \qed \end{example}
\begin{example}Let $q\p\sim 1$ in $\mathcal{C}_{p,q}^{M}$ and
  $\mathcal{C}_{p,qq\p}^{M}$.
  If $\mathcal{C}_{p,q}^{M}$ contains single equivalence class,
  $H_{p,qq\p}^{M}\cong \{1\}$.
  If $\mathcal{C}_{p,qq\p}^{M}$ contains single equivalence class,
  $H_{p,q}^{M}\cong \{1\}$. \qed\end{example}
\begin{example}$H_{7,16}^{M}$ and $H_{7,2}^{M}$ are isomorphic to quotient groups of
   each other as $8\in \mathrm{Ker\: }\phi_{7,2}\cap \mathrm{Ker\:
   }\phi_{7,16}$.\\
   $H_{5,2}^{M}$ and $H_{5,12}^{M}$ are isomorphic to quotient groups of
   each other as $6\in \mathrm{Ker\: }\phi_{5,2}\cap \mathrm{Ker\:
   }\phi_{5,12}$.\\
   $H_{3,2}^{M}$ and $H_{3,8}^{M}$ are isomorphic to quotient groups of
   each other as $4\in \mathrm{Ker\: }\phi_{3,2}\cap \mathrm{Ker\:
   }\phi_{3,8}$. \qed
 \end{example}
 
\begin{example}Let $q\p \in \mathrm{Ker\: }\phi_{p,q} \cap \mathrm{Ker\: }\phi_{p,qq\p}$.
  If $\mathcal{C}_{p,qq\p}^{M}$ contains one equivalence
  class, then $a\nsim 1$ for $a\in M$ in $\mathcal{C}_{p,q}^{M}$ only if $q\p|a$.\\
  If $\mathcal{C}_{p,q}^{M}$ contains one equivalence
  class, then $a\nsim 1$ for $a\in M$ in $\mathcal{C}_{p,qq\p}^{M}$ only if $q\p|a$.
\end{example}
\proof
Let $\mathcal{C}_{p,qq\p}^{M}$ contain a single equivalence
  class. Then $H^{M}_{p,qq\p}\cong \{1\}$ and $\overline{H}^{M}_{p,qq\p}\cong F$.
We have $$\nicefrac{\overline{H}^{M}_{p,qq\p}}{\langle q\p=1\rangle}\cong
\nicefrac{\overline{H}^{M}_{p,q\p}}{\langle q\p=1\rangle}\cong
\nicefrac{F}{\langle q\p=1\rangle}$$ by Property \ref{qp-factorization}.
Therefore all members of $M$ that are not equivalent to $1$ are not
co-prime to $q\p$. Similar argument is applied for the case when $\mathcal{C}_{p,q}^{M}$ contains a single equivalence
  class.
 \qed
\bibliographystyle{plain}

\end{document}